\documentstyle{amsppt}
\topmatter
\author Dale E. Alspach \endauthor
\thanks Research supported in part by National Science Foundation
grant DMS-8902327\endthanks
\footnote""{This paper is in final form and no version of it will be\newline
submitted for publication elsewhere.}
\address Department of Mathematics, Oklahoma State University,
\newline Stillwater, OK 74078-0613
\endaddress
\email alspach\@hilbert.math.okstate.edu \endemail
\title
A $\ell_1$-Predual which is not Isometric to a Quotient of $C(\alpha )$
\endtitle
\abstract
About twenty years ago Johnson and Zippin showed that every
separable $\text{L}_1(\mu)$-predual was isometric to a quotient of
$C(\Delta
)$, where
$\Delta$ is the Cantor set. In this note we will show that the
natural analogue of the theorem for $\ell_1$-preduals does not
hold. We will show that there are many $\ell_1$--preduals which
are not isometric to a quotient of any $C(K)$-space with $K$
a countable compact metric space.
We also prove some general results about the relationship
between $\ell_1$-preduals and quotients of $C(K)$-spaces with $K$ a
countable compact metric space. 
\endabstract
\endtopmatter
\document
About 20 years ago Johnson and Zippin \cite{J-Z} proved the following theorem.

\proclaim{Theorem} Suppose that $X$ is a separable
$\text{L}_1(\mu)$-predual, then there is subspace $Y$ of $C(\Delta 
)$, where 
$\Delta$ is the Cantor set, such that $X$ is isometric to
$C(\Delta)/Y$.
\endproclaim

Because the space $X$ might be $C(\Delta)$, $C(\Delta)$ is the
smallest $\text{L}_1(\mu)$-predual that one could use for such a result. If we
consider the class of $\ell_1$-preduals then it is conceivable
that a smaller space might be sufficient although the space
would need to depend on some measurement of the size of the
$\ell_1$-predual. A natural class of spaces to consider is the
spaces $C(\alpha)$ where $\alpha$ is a countable ordinal. (This
is the same as the class of $C(K)$-spaces with $K$ a countable
compact Hausdorff space by a classical result of Mazurkiewicz and
Sierpi\'nski, \cite{M-S}.) Thus one can consider the following question.

\proclaim{Question} If $X$ is an $\ell_1$--predual is there a
countable ordinal $\alpha$ such that $X$ is isometric (isomorphic)
to a quotient of $C(\alpha)$?
\endproclaim

We will show that the isometric question has a negative answer and prove
some technical results which are useful for deciding whether an
$\ell_1$--predual is isomorphic to a quotient of $C(\alpha)$. The
isomorphic question remains open and at the end of the paper we
discuss some variants of the isomorphic problem. Note also that
we consider only $X$ for which $X^*$ is isometric to $\ell_1$ because
even the Johnson--Zippin result is false for isomorphic
$\ell_1$-preduals, \cite{B-D}.

Throughout this paper $\ell_1$-predual will mean a Banach space
with dual isometric to $\ell_1$. If $\alpha$ is an ordinal
$C(\alpha)$ will denote the space of continuous functions on the
ordinals less than or equal to $\alpha$ with the order topology.
If $A$ is a subset of a Banach space, $[A]$ is the norm closed
linear span of $A$. Notation and standard results from Banach
space theory may be
found in the books of Lindenstrauss and Tzafriri, \cite{L-T}, Dunford
and Schwartz, \cite{D-S}, and Diestel, \cite{D}.

We
begin with some technical results about basic sequences
equivalent to the usual $\ell_1$ basis which are contained
in dual spaces.

\proclaim{Lemma 1} Let $Z$ be a Banach space with separable dual
and let $Y$ be a subspace of $Z^*$ which is isomorphic to $\ell_1$
with normalized $\ell_1$-basis $(y_n)$. If \ $\overline{\{y_n\}}^{\text{w}^*}
\subset Y$ then $Y$ is $\text{w}^*$-closed in $Z^*$.
\endproclaim

\demo{Proof} By the Krein-Smulian Theorem, \cite{D-S, Theorem V.5.7},
it is sufficient to
show that $B_Y=\{y\in Y: \|y\|\leq 1\}$ is
$\text{w}^*$-closed. Let $(z^*_n)$ be a sequence in $B_Y$
with $\text{w}^*$-limit $z^*$.  We will show that $z\in Y.$ Let
$K$=$\overline{\{y_n\}}^{\text{w}^*}$ and let $\lambda$ 
$\geq1$ satisfy $\lambda \|\sum a_n y_n\| \geq \sum|a_n|$
for all finite sequences $(a_n)$. We know that
$\overline{\text{co} \pm \lambda
K}^{\text{w}^*}\supset\overline{\text{co}\{ \pm \lambda y_n\}}^{\text{w}^*}
\supset B_Y.$

Because $z^* \in \overline{B}^{\text{w}^*}_Y \subset
\overline{\text{co} \pm \lambda
K}^{\text{w}^*},$ by the Choquet, \cite{D, p 154}, and Milman Theorems,
\cite{D-S, Lemma
V.8.5},  $z^*$ is
represented by a probability measure $\mu$ supported on
$$\text{Ext }\overline{\text{co} \pm \lambda 
K}^{\text{w}^*}\subset \pm\lambda K.$$
If $z^* \notin Y$ then by the
Hahn-Banach Theorem there exists $z^{**}\in Z^{**}$ such that
$z^{**}_{|Y}=0$, $\|z^{**}\|=1$, and $z^{**}(z^*)>0$. $z^{**}$ is
the $\text{w}^*$-limit of a sequence in $Z$, $(z_j)$, with
$\|z_j\|\leq 1$ for all $j$. Thus
$$ z^{**}(z^*)=\lim_j z^*(z_j)=\lim_j\int_{\pm \lambda K} z_j
d\mu.$$
Since $K\subset Y$, $z^{**}(k)=0$ for all $k\in K$. Hence
$z_j(\lambda k)\rightarrow 0$ for all $k\in K$. Therefore $(z_j)$
is a uniformly bounded sequence converging pointwise to 0 on $\pm
\lambda K$ and it follows from the Bounded Convergence Theorem
that $z^{**}(z^*)=0.$ \qed \enddemo

The next lemma tells us that the $\text{w}^*$-closure of the $\ell_1$-basis is
the only thing that is important.

\proclaim{Lemma 2} Suppose that $X$ and $Y$ are separable
Banach spaces and
that $(x_n^*)$ and $(y_n^*)$ are normalized sequences in $X^*$
and $Y^*$, respectively,
which are equivalent to the usual unit vector basis of
$\ell_1$ and for which $\overline{[x_n^*]}^{\text{w}^*}=[x_n^*]$
and $\overline{[y_n^*]}^{\text{w}^*}=[y_n^*]$.
Suppose that
the basis to basis map $\phi$ of $[x_n^*]$ onto
$[y_n^*]$, i.e., $\phi(\sum a_n x_n^*)=\sum a_n y_n^*$,
 is a $\text{w}^*$-homeomorphism
of the $\text{w}^*$-closure of $\{x_n^*\}$ onto the
$\text{w}^*$-closure of
$\{y_n^*\}$.
Then $\phi$ is a $\text{w}^*$-continuous isomorphism of $[x_n^*]$
onto $[y_n^*]$.
\endproclaim

\demo{Proof} Only the $\text{w}^* $-continuity
needs to be proved. By passing to
quotients we may assume that $(x_n^*)$ and $(y_n^*)$ are bases
for the duals
of $X$ and $Y$, respectively. If we identify X and Y with the $\text{w}^*
$-continuous affine symmetric functions on $B_{X^*}$ and
$B_{Y^*}$, respectively,
we need only show that the linear extension of a $\text{w}^*
$-continuous
linear function on $K(X)=\lambda^{-1}\overline{\{x_n^*\}}^{\text{w}^*} \cup
\{0\}$, where $(x_n^*)$ is $\lambda $ equivalent to the usual
basis of $\ell_1$,
is $\text{w}^*
$-continuous on $B_{X^*}$.
Let $\Phi$ be the map defined by $(\Phi f)(x^*)=f(\phi(x^*))$
for $x^*\in K(X)$ and let $J_X$ and $J_Y$ be the evaluation maps from
$X$ into $C(K(X)$ and from $Y$ into $C(K(Y)$, respectively. ($K(Y)$
is defined analogously to $K(X)$.)
We have the following diagram.
$$
\CD Y @. X\\
@V{J_Y}VV @VV{J_X}V \\
C(K(Y)) @>\Phi>> C(K(X))
\endCD
$$

In order to invert $J_X$
we need to show that any $x^{**}\in C(K(X))\cap X^{**}$ is actually
in $X$. Any such $x^{**}$ is the $\text{w}^*$--limit of a sequence
$(x_n)\in \|x^{**}\|B_X$. Because $x^{**}$ is continuous
on $K(X)$, $(x_n)$ converges to $x^{**}$ weakly in $C(K(X))$. Hence
$x^{**}_{|K(X)}\in \overline{\text{co } \{{x_n}_{|K(X)}\}}^{\|\cdot\|}$ and
thus $x^{**}\in X$ because $K(X)$ is norming.
\qed \enddemo

If we reformulate these results for the special case where $X$ is
a $\ell_1$-predual and $Y$ is some $C(\alpha)$, we get

\proclaim{Proposition 3} Suppose that $X$ is an $\ell_1$-predual
and let $(x_n^*)$ be
the $\ell_1$-basis of the dual. Then $X$ is isometric to a quotient of
$C(\alpha)$, for some $\alpha<\omega_1$,
if and only if there is a sequence $(\mu_n)$ 
of norm 1,
mutually singular measures on $\alpha$ such that
$K=\overline{\{\mu_n^*\}}^{\text{w}^*}\subset [\mu_n]$ and the map $\phi$
defined by $\phi(x_n^*)=\mu_n$ extends to a $\text{w}^*$-continuous 
map from
$\overline{\{x_n^*\}}^{\text{w}^*}$ onto
$K$.\endproclaim

\demo{Proof} It is well known that a sequence of measures in the dual of a
$C(K)$-space is isometric to a $\ell_1$-basis if and only if the
measures are norm one and mutually singular. The two lemmas
complete the proof.
\qed\enddemo

To our knowledge this proposition has not appeared previously but
for a long time it was used on an intuitive level by the author.
This can be seen in the construction of the example in \cite{A} and
in the treatment of the Benyamini-Lindenstrauss example \cite{B-L}
given in \cite{A-B}. To see the utility of the proposition consider
the following standard example of an $\ell_1$-predual which is
not isometric to a $C(K)$-space.

\example{Example 1}
Let $(t_n)$ be a sequence of points in the open unit square in $\Bbb
R^2$ such that $K=\overline{\{t_n\}}=\{t_n\}\cup\{(t,0):0\leq t\leq
1\}.$ Let $X=\{f\in C(K):f(t,0)=tf((0,0))+(1-t)f((1,0))\}$. Clearly
$X$ is isomorphic to $c_0$ and the $\ell_1$-basis for the dual is
the evaluations at the $t_n$'s and at (0,0) and (1,0).
We claim that $X$ is isometric to
a quotient of $C(\omega\cdot 2)$. Indeed, let $t_{n,1}$ denote
the first coordinate of $t_n$ for all $n$.  Define
$\phi(\delta_{t_n})=t_{n,1}\delta_n+(1-t_{n,1})\delta_{\omega+n}$,
$\phi(\delta_{(0,0)})=\delta_\omega $,
and $\phi(\delta_{(1,0)})=\delta_{\omega\cdot 2} $. 
It is easy to verify that $\phi$ extends to be a
$\text{w}^*$-continuous map from the closure of the basis of
$X^*$ onto the closure of the image. Clearly the span of the
range of $\phi$ contains the $\text{w}^*$-closure of the range of
$\phi$. Thus $X$ is a quotient of $C(\omega\cdot 2)$ by the
proposition.
\endexample

Next we will note some simple facts about basic sequences
in the dual of a $C(K)$-space which are 1-equivalent to the usual
basis of 
$\ell_1$.

\proclaim{Lemma 4} If $(\mu_n)$ is a sequence of mutually
singular (non-zero) measures in $C(K)^*$, $K$ compact metric, which
converges $\text{w}^*$ to $\mu$ then $$\text{supp }\mu\subset
\{k\in K:\text{there exist } k_n\in \text{supp }\mu_n \text{with
} \lim k_n = k\}.$$
\endproclaim

If the measures in Lemma 4 are all atomic then we may
replace the support of the measures by the set of points where
the measure is non-zero and the result remains true. The proof of
the lemma is
straight-forward so we leave it to the reader.

\proclaim{Corollary 5} If $K$ is a countable compact metric space
and $(\mu_n)$ is a sequence of mutually
singular measures in $C(K)^*$ with $\text{w}^*$-limit $\mu$,
then $\text{supp } \mu \subset
K^{(1)}$, the first derived set of $K$.
\endproclaim

\demo{Proof} This follows immediately from Lemma 4. \qed \enddemo

Now let $X$ be a $\ell_1$-predual and let $(e_n)$ be the
$\ell_1$-basis for the dual. We want to define a system of
derived sets $N^{(\alpha)}=N^{(\alpha)}((e_n)_{n\in \Bbb N})$
 of $\Bbb N$ that has properties similar to those
of the Szlenk  sets
for the
$\text{w}^*$-closure $K$ of the sequence $(e_n)$.
(See \cite{S} for the definition and properties of the Szlenk sets.)

\definition{Definition} Let $N^{(0)}=\Bbb N$ and if $N^{(\alpha)}$ has been
defined, let 
$$\multline
N^{(\alpha+1)}
=\{n\in N^{(\alpha)}:\text{there exists
an infinite }
M\subset N^{(\alpha)}\\
 \text{ such that } (\text{w}^* \lim_{j \in M} e_j)
(n)\neq 0\}.
\endmultline$$
 If $\alpha$ is a limit ordinal define
$N^{(\alpha)}=\cap_{\beta<\alpha} N^{(\beta)}.$
\enddefinition

\proclaim{Proposition 6} Let $K$ be a countable compact metric
space. If $(\mu_n)$ is a sequence of norm one
measures in $C(K)^*$, $K$ a countable compact metric space,
which are 1-equivalent to the unit vector
basis of $\ell_1$ and $[\mu_n]$ is $\text{w}^*$-closed, then
$\{N^{(\alpha)}((\mu_n)_{n\in\Bbb N}):\alpha< \alpha_0\}$
is a strictly decreasing family of subsets  of $\Bbb N$,
where $\alpha_0$ is the smallest ordinal $\alpha$
such that $N^{(\alpha)}=\emptyset.$
\endproclaim

\demo{Proof} By Corollary 5 we have that if $(\mu_j)_{j\in M}$
converges to a non-zero measure $\mu=\sum a_k \mu_k$, then $\text{supp }\mu
\subset K^{(1)}$. Therefore if $a_n\neq 0$ then $n\in N^{(1)}$,
and $\text{supp }\mu\subset K^{(1)}$ implies that 
$\text{supp }\mu_n
\subset K^{(1)}$ because the measures $\mu_k$ are mutually
singular. Hence for all $n\in  N^{(1)}$, $\text{supp }\mu_n
\subset K^{(1)}$. A simple induction argument
completes the proof. \qed \enddemo

Using this proposition it is very easy to find $\ell_1$-preduals
which are not isometric to a quotient of any $C(\alpha)$.

\example{Example 2} Let $\dsize X=\{(a_n)\subset \Bbb R:\lim_{j
\rightarrow \infty} a_j=\sum_{n\in \Bbb N}\frac{a_n}{2^n}\}.$ Note
that $X$ is
a codimension 1 subspace of $c=C(\omega).$ Thus $X$ is
isomorphic to $c_0$. We claim that the evaluation functionals defined by
$e_n((a_j))=a_n$ are a basis for the dual 1-equivalent to the usual
$\ell_1$-basis. To do
this for each $k\in \Bbb N$ we will exhibit  a copy
of $\ell_\infty^k$ which norms the
span of the first $k$ of these functionals. Fix $K>k$ and let
$$f_j(n)=\cases 1, &\text{if }n=j\\0, &\text{if }n < K \text{ and }n\neq j\\
-2^{-j},&\text{if }n=K\\2^{-j},&\text{if }n>K.\endcases$$
for $j=1,2,\dots,k.$ Obviously $\lim_{n\rightarrow \infty}
 f_j(n)=2^{-j}=\sum_{n\in
\Bbb N} 2^{-n} f_j(n).$ Thus $(f_j)_{j\leq k}$ is a sequence in $X$ 1-equivalent to the
unit vectors in $\ell_\infty^k$  and biorthogonal to the first $k$
evaluation functionals.

We claim that $X$ is not isometric to a quotient of $C(\alpha)$.
Indeed, the $\text{w}^*$-limit in $X^*$ of $(e_j)$ is $\dsize 
\sum_{n\in \Bbb N}
\frac{e_n}{2^n}.$ Hence $N^{(\alpha)}((e_n)_{n\in \Bbb N})=\Bbb N$
for all $\alpha$. By Propositions 3 and 6,  $X^*$ is not
isometric to a $\text{w}^*$-closed subspace of $C(K)^*$ for any
countable compact metric space $K$, consequently, $X$ is not
isometric to a quotient of any $C(\alpha).$

The point here is that in this case one really needs the Cantor
set or some other uncountable compact metric space to support the
measures. If we relax the requirement from isometric to a
quotient to
$(1+\epsilon)$-isomorphic to a quotient of $C(\alpha)$, then we
can embed the basis of $X^*$ in $C(\omega^2)^*$.

Fix $N$ and define $\phi(e_n)=\delta_{\omega\cdot n}$ for $n\leq
N$ and
$$\align
\phi(e_n)&=(1-\frac{1}{2^{n-1}})^{-1}[\sum_{j\leq N}
\frac{\delta_{\omega\cdot (j-1)+n}}{2^j}
+ \sum_{N+1\leq j<n}
\frac{\phi(e_j)}{2^j}]\\
&=(1-\frac{1}{2^{n-1}})^{-1}\sum_{1\leq j<n}
\frac{\phi(e_j)}{
2^j}
\endalign$$
for $n> N.$
It follows from Lemma
2 that $\phi$ can be extended to a $\text{w}^*$-isomor\-phism of
$X^*$ into $C(\omega^2 )$. \endexample

The example provides a trivial way of producing more such
examples. Indeed, take any $\ell_1$-predual
and construct an isomorph by adding the example as a
$\ell_\infty$ direct summand to get another $\ell_1$-predual
which is not isometric to a quotient of any $C(\alpha)$. There
are probably much more interesting and subtle examples that can
be constructed by using the behavior of the $\text{w}^*$-topology
in the example.

As we noted earlier it is unknown whether every $\ell_1$-predual
is isomorphic to a quotient of $C(\alpha)$. In addition there are
several variants of this question which are open. Below
$X$ is a $\ell_1$-predual.

\example{Question 1}   Given $\epsilon>0$ does there exist
$\alpha<\omega_1$ and a $\text{w}^*$-continuous into isomorphism
$S$ of $X^*$ into $C(\alpha)^*$ with $\|S\|
\|S^{-1}_{|S(X^*)}\|<1+\epsilon$? In other words is $X$
$(1+\epsilon)$-isomorphic to a quotient of $C(\alpha)$?
\endexample

\example{Question 2} Is there an $\alpha<\omega_1$ and a
quotient Y of $C(\alpha)$ such that $Y^*$ is isometric to
$\ell_1$ and $X$ is isomorphic to $Y$?
\endexample

A positive answer to Question 1 or 2 would of course answer the
isomorphic question. A positive answer to Question 1 would imply
that there is an $\alpha$ that works for all $\epsilon$
simultaneously, that is, there is an $\alpha$ so that $X$ is
almost isometric to a quotient of $C(\alpha)$. Such a result should
provide very useful information of the structure of
$\ell_1$-preduals. Let us also note that very little is known
about quotients of $C(\alpha)$. In particular it is unknown whether
they must be $c_0$-saturated. ($Y$ is $c_0$-saturated if every
subspace contains $c_0$.)

\enddocument
\bye